\documentstyle[11pt]{article}
\def\Bbb R{{\rm \bf R}}
\def\proclaim#1{\vskip2mm{\bf #1}\em}
\def\endproclaim{\em \vskip2mm}
\def\tag#1{\eqno(#1)}
\def\gathered{\begin{array}{c}}
\def\endgathered{\end{array}}
\def\text{\mbox}

\begin{document}

\title {On finding an obstacle embedded in the rough background medium via the enclosure method in the time domain}
\author{Masaru IKEHATA\footnote{
Laboratory of Mathematics,
Institute of Engineering,
Hiroshima University, Higashihiroshima 739-8527, JAPAN}}
%\date{}
\maketitle

\begin{abstract}
A mathematical method for through-wall imaging via wave phenomena in the time domain
is introduced.  The method makes use of a single reflected wave over a finite time interval
and gives us a criterion whether a penetrable obstacle exists or not
in a general rough background medium.
Moreover, if the obstacle exists, the lower and upper estimates of the distance between the obstacle
and the center point of the support of the initial data are given.
As an evidence of the potential of the method two applications are also given.

%{\bf FinalRevised12SubmittedIkehata.ImagingBehindWall.tex}

\noindent
AMS: 35R30

\end{abstract}

%\tableofcontents

\section{Introduction}

Assume that there is a large wall between an observer and an
unknown object. The wall is {\it opaque} and thus the observer can
not see the object directly.  How can the observer find the
object? Consider the case when the wall is electromagnetically
{\it penetrable}.  Then there is no doubt that everyone thinks
about the use of electromagnetic wave.  Generate the
electromagnetic wave at the place where the observer is.
And observe the reflected wave from the object at the same place over
a finite time interval. The observed wave should include
information about the object. How can one extract the information
from the observed wave? The purpose of  this paper is to develop a
mathematical method for {\it through-wall imaging} by using the
governing equation of the wave from the beginning to end.
Originally the governing equation should be the Maxwell system. In
this paper, as a first step we consider two single equations for
{\it scalar waves}.

\subsection{Finding discontinuity in a medium with a rough refractive index}

Let us formulate the first problem.
Let $n=1,2,3$.
Let $\alpha\in L^{\infty}(\Bbb R^n)$ and satisfy $\text{ess.inf}_{x\in\Bbb R^n}\alpha(x)>0$.
Let $0<T<\infty$.  Given $f\in\,L^2(\Bbb R^n)$, let
$u=u_f(x,t)$ be the weak solution of
$$
\left\{
\begin{array}{ll}
\displaystyle
\alpha(x)\partial_t^2u-\triangle u=0 & \text{in}\,\Bbb R^n\times\,]0,\,T[,
\\
\\
\displaystyle
u(x,0)=0 & \text{in}\,\Bbb R^n,\\
\\
\displaystyle
\partial_tu(x,0)=f(x)  & \text{in}\,\Bbb R^n.
\end{array}
\right.
\tag {1.1}
$$
The notion of the weak solution is taken from \cite{DL}.  See Subsection 2.1 in this paper.

We assume that $\alpha$ takes the form
$$\displaystyle
\alpha(x)=\left\{
\begin{array}{ll}
\displaystyle
\alpha_0(x), & \mbox{if $x\in \Bbb R^n\setminus D$},\\
\\
\displaystyle
\alpha_0(x)+h(x), & \mbox{if $x\in D$},
\end{array}
\right.
\tag {1.2}
$$
where $D\subset\Bbb R^n$ is a bounded open subset with Lipschitz boundary.
The function $\alpha_0$ belongs to $L^{\infty}(\Bbb R^n)$
and satisfies $m_0^2\le\alpha_0(x)\le M_0^2$ a.e.$x\in\Bbb R^n$ for positive constants $m_0$ and $M_0$;
$h$ belongs to $L^{\infty}(D)$ and satisfies one of (A.I) and (A.II) listed below:

(A.I) $\exists C>0\,\,h(x)\ge C$ a.e. $x\in D$;

(A.II) $\exists C>0\,\, -h(x)\ge C$ a.e. $x\in D$.

$D$ is a model of an unknown {\it penetrable obstacle} with material parameter $\alpha_0+h$
embedded in the background medium with material parameter $\alpha_0$.  The distribution of $\alpha_0$
models various penetrable walls in the space and $D$ is something hidden in the walls or a space
between the walls and various penetrable obstacles.

Let $B$ be an open ball satisfying $\overline B\cap\overline
D=\emptyset$. Let $f\in L^2(\Bbb R^n)$ satisfy $f(x)=0$ a.e.
$x\in\Bbb R^n\setminus B$ and $\text{ess.inf}_{x\in B}\,f(x)>0$.
Generate $u_f$ by the initial data $f$. In this paper, we consider
the following inverse problems under the assumption that
$\alpha_0$ is {\it known} and that both $D$ and $h$ are {\it
unknown}.

\noindent
{\bf Problem 1.}
Find a criterion whether $D=\emptyset$ or not in terms of only $u_f$ on $B$ over
time interval $]0,\,T[$.

\noindent
{\bf Problem 2.}
Assume that $D\not=\emptyset$.
Extract information about $D$ from $u_f$ on $B$ over time interval $]0,\,T[$ for a fixed large $T$.

Note that the correspondence $(D,h)\longmapsto u\vert_{B\times\,]0,\,T[}$ is nonlinear
and thus both problems become nonlinear problems.
The existence of variation of $\alpha_0$ outside $D$ plays a role of the wall in which a wave can propagate.
This paper aims at developing a {\it through-the wall} imaging method via the governing equation
on the wave used.

There are a lot of studies from the engineering side for through-the-wall imaging using electromagnetic waves.
See \cite{B, AA} for a survey on through-wall imaging and various engineering approaches.
Roughly speaking, one can say that those approaches introduce various {\it processing} of the reflected {\it signal} from the wall
and targets behind the wall. For example, in \cite{AA} the {\it compressive sensing} incorporating invariance of the antenna location due to
the geometry of the assumed wall has been applied to a wall clutter mitigation technique for the signal.
In \cite{AZA} under the assumption that the wall is a single uniform one,
an approach which employs an imaging function
incorporating {\it geomtrical optics} (Snell's law) for the wave propagation through the wall
has been considered and tested numerically. 
In \cite{Z} an algorithm to find a moving human in a simple situation using the time-of-flight approach is introduced.
In \cite{CC} the idea of the {\it linear sampling method} in the {\it frequency domain} has been applied to through-wall imaging
and tested numerically in two dimensions under the assumption that the wall is infinitely long in one direction. 
They employ the concrete form of the Green's function for the wave propagation through the background medium and thus, in this sense,
their approach should belong to a class of PDE approaches.

In this paper, we employ the Enclosure Method introduced by the
author himself in \cite{I1,IE} as a guiding principle for
attacking the problem mentioned above. It is a direct method which
connects the unknown discontinuity and the observation data
through the so-called indicator function. In \cite{I4} the author
initiated to apply the method to inverse obstacle problems whose
governing equations are given by the heat and wave equations in
one-space dimension. Now we have several applications of the
method to inverse obstacle scattering problems in three-space
dimensions whose governing equations are given by the wave
equations \cite{IE0, IW2, IW3, IWCOE}. See \cite{IR} for a review
of the results in those papers. However, imaging an obstacle
through-the wall is still out of the range of the results in those
papers.

Now let us describe the result.

Let $\tau>0$ and define
$$\displaystyle
w(x,\tau)=\int_0^T e^{-\tau t}u(x,t)dt\,\,\text{in}\,\Bbb R^n.
\tag {1.3}
$$
Let $v=v(x,\tau)\in H^1(\Bbb R^n)$ be the weak solution of
$$\displaystyle
\triangle v-\alpha_0\tau^2 v+\alpha_0 f=0\,\,\text{in}\,\Bbb R^n.
\tag {1.4}
$$
Define
$$\displaystyle
I_f(\tau,T)
=\int_B\alpha_0f(w-v)dx.
$$

The following result gives us solutions to Problems 1 and 2.

\proclaim{\noindent Theorem 1.1.}
We have:

(i) if $D=\emptyset$, then for all $T>0$ it holds that
$$\displaystyle
\lim_{\tau\longrightarrow\infty}e^{\tau T}I_f(\tau,T)=0;
$$

(ii) if $D\not=\emptyset$ and $h$ satisfies (A.I), then for all $T>2M_0\text{dist}(D,B)$ it holds that
$$\displaystyle
\lim_{\tau\longrightarrow\infty}e^{\tau T}I_f(\tau,T)=-\infty;
$$

(iii)  if $D\not=\emptyset$ and $h$ satisfies (A.II), then for all $T>2M_0\text{dist}(D,B)$ it holds that
$$\displaystyle
\lim_{\tau\longrightarrow\infty}e^{\tau T}I_f(\tau,T)=\infty.
$$

Moreover, in case of both (ii) and (iii) we have, for all $T>2M_0\text{dist}(D,B)$
$$\displaystyle
\liminf_{\tau\longrightarrow\infty}
\frac{1}{2\tau}
\log\left\vert I_f(\tau,T)\right\vert
\ge -M_0\text{dist}\,(D,B)
\tag {1.5}
$$
and
$$\displaystyle
\limsup_{\tau\longrightarrow\infty}
\frac{1}{2\tau}
\log\left\vert I_f(\tau,T)\right\vert
\le -m_0\text{dist}\,(D,B).
\tag {1.6}
$$

\endproclaim

Note that $\alpha_0$ and $h$ are just {\it essentially bounded} on $\Bbb R^n$ and $D$, respectively.
In particular, $\alpha_0$ can be a model for various background media such as multilayered media with complicated interfaces or unions of various domains with different refractive indexes.
It will be {\it impossible} to apply any approach based on the idea of {\it geometrical optics} to this case.
See \cite{H, MT, P1, P2} for such approach in the case when the {\it scattering kernel} which is the observation data
in the Lax-Phillips scattering theory is given under the assumptions that $\alpha_0(x)=1$ a.e. $x\in\Bbb R^n$, $\partial D$ is smooth
and $h\in C^{\infty}(\overline D)$.

Let $p$ be the center of $B$ and $\eta$ the radius. 
Define $d_{\partial D}(p)=\inf_{x\in\partial D}\vert x-p\vert$.
Since $\text{dist}\,(D,B)=d_{\partial D}(p)-\eta$, estimates (1.5) and
(1.6) give us an upper and lower estimate of $d_{\partial D}(p)$
provided $m_0$ and $M_0$ are known and $T$ is sufficiently large.
Thus we obtain the upper and lower estimation of the sphere $\vert
x-p\vert=d_{\partial D}(p)$ whose exterior {\it encloses} the
object. Estimates (1.5) and (1.6) suggest a {\it new direction} of
the Enclosure Method in the case when the background medium is
{\it inhomogeneous} and quite complicated: give up to find a
precise quantity in the observation data which is related to the
exact location of unknown obstacles; instead give lower and upper
estimates {\it rigorously} like (1.5) and (1.6) for
$\text{dist}\,(D,B)$.

Some further remarks are in order.

$\bullet$  In Theorem 1.1 it suffices to know $v$ on $B$ not whole
$v$. However, without knowing $\alpha_0$ outside $B$ it is
impossible to compute $v$ on $B$ in advance. In the last section
we suggest an experimental computation procedure of $v$ on $B$
without detailed knowledge of $\alpha_0$ outside $B$.  It seems
that this will be useful for the daily security of the interior of
the room, house and building which have complicated environment.
For the purpose Theorem 1.1 will be suitable since $\alpha_0$ is
just essentially bounded and we never assume any other regularity.

$\bullet$  If $\alpha_0(x)=1$ a.e. $x\in\Bbb R^n$, then one can choose $M_0=m_0=1$ and thus (1.5) and (1.6) imply that
$$\displaystyle
\lim_{\tau\longrightarrow\infty}
\frac{1}{2\tau}
\log\left\vert I_f(\tau,T)\right\vert
=-\text{dist}\,(D,B).
$$
This coincides with a result (1.20) in \cite{IW2}.

In the next subsection we apply the idea developed here to more special case and show that (1.5) and (1.6) can be replaced with a single formula.

\subsection{Finding discontinuity in a dissipative medium}

Let $n=1,2,3$.
Let $q\in L^{\infty}(\Bbb R^n)$ satisfy $q(x)\ge 0$ a.e.$x\in\Bbb R^n$.
Let $0<T<\infty$.
Given $f\in\,L^2(\Bbb R^n)$, let
$u=u_f(x,t)$ be the weak solution of
$$
\left\{
\begin{array}{ll}
\displaystyle
\partial_t^2u-\triangle u+q(x)\partial_tu=0 & \text{in}\,\Bbb R^n\times\,]0,\,T[,
\\
\\
\displaystyle
u(x,0)=0 & \text{in}\,\Bbb R^n,\\
\\
\displaystyle
\partial_tu(x,0)=f(x) & \text{in}\,\Bbb R^n.
\end{array}
\right.
\tag {1.7}
$$
We assume that $q$ takes the form
$$\displaystyle
q(x)=\left\{
\begin{array}{ll}
\displaystyle
q_0(x), & \,\mbox{if $x\in\Bbb R^n\setminus D$}\\
\\
\displaystyle
q_0(x)+h(x), & \,\mbox{if $x\in D$},
\end{array}
\right.
\tag {1.8}
$$
where $D$ is a bounded open set of $\Bbb R^n$ with Lipschitz boundary.
The function $q_0$ belongs to $L^{\infty}(\Bbb R^n)$ and satisfies $q_0(x)\ge 0$ a.e. $x\in\Bbb R^n$.

The function $h$ on $D$ in (1.8) has to belong to $L^{\infty}(D)$.
We assume that $h$ satisfies one of (A.I) and (A.II).

Let $\tau>0$.  Define $w$ by the formula (1.3) in which $u$ is replaced with the solution of (1.7).

Let $v\in H^1(\Bbb R^n)$ be the weak solution of
$$\displaystyle
(\triangle-\tau^2-\tau q_0)v+f=0\,\,\text{in}\,\Bbb R^n.
\tag {1.9}
$$
Let $B$ and $f$ be the same as those of Theorem 1.1.
Define
$$\displaystyle
J_f(\tau,T)
=\int_Bf(w-v)dx.
$$
The following result is new and suggests that, even in a rough dissipative medium the solution
of (1.7) carries information  about $D$ clearly like the wave equation outside $D$.

\proclaim{\noindent Theorem 1.2.}  
We have:

(i) if $D=\emptyset$, then for all $T>0$ it holds that
$$\displaystyle
\lim_{\tau\longrightarrow\infty}e^{\tau T}J_f(\tau,T)=0;
$$

(ii) if $D\not=\emptyset$ and $h$ satisfies (A.I), then for all $T>2\text{dist}(D,B)$ it holds that
$$\displaystyle
\lim_{\tau\longrightarrow\infty}e^{\tau T}J_f(\tau,T)=-\infty;
$$

(iii)  if $D\not=\emptyset$ and $h$ satisfies (A.II), then for all $T>2\text{dist}(D,B)$ it holds that
$$\displaystyle
\lim_{\tau\longrightarrow\infty}e^{\tau T}J_f(\tau,T)=\infty.
$$

Moreover, in case of both (ii) and (iii) we have, for all $T>2\text{dist}(D,B)$
$$\displaystyle
\lim_{\tau\longrightarrow\infty}
\frac{1}{2\tau}
\log\left\vert J_f(\tau,T)\right\vert
=-\text{dist}\,(D,B).
\tag {1.10}
$$

\endproclaim

Roughly speaking, we see that (1.9) as $\tau\longrightarrow\infty$ corresponds to (1.4) with $\alpha_0=1$ and thus
formula (1.10) is reasonable.

 Comparing Theorem 1.2 with the previous results in \cite{IW2, IW3} for the wave equation outside $D$,
 we see
 (ii) and (iii) suggest us that that assumptions (A.I) (stronger dissipation) and (A.II) (weaker dissipation)
 correspond to the Dirichlet and Neumann boundary conditions imposed on $\partial D$, respectively.

A brief outline of this paper is as follows.
In Section 2 first we specify the meaning of the weak solution used in the formulation of the problems.
Second we establish two basic integral identities for $w$ given by (1.3).
Theorem 1.1 is proved in Section 3.
The proof starts with deriving the lower and upper estimates for the indicator function:
$\tau\longmapsto I_f(\tau,T)$
as $\tau\longrightarrow\infty$ from the basic identities in
Section 2. Next, by virtue of the governing equation (1.4) we see
that $v$ has point-wise explicit lower and upper estimates.
Combining those, we obtain all the conclusions stated in Theorem
1.1. The proof of Theorem 1.2 is given in Section 4 which is a
combination of a reduction and similar argument done in the proof
of Theorem 1.1. In Section 5 we give some remarks concerned with a
``practical use'' of Theorem 1.1 and present a sharp form of
Theorem 1.1 in the case when the space dimension is $1$ and
$\alpha_0$ has a special but important form.

\section{Preliminaries}

\subsection{On the solution class}

In this subsection we specify the meaning of the weak solutions of (1.1) and (1.7) at the same time.
It is based on Theorem 1 given on p.558 in \cite{DL} which we have already used for the formulation
of the weak solution of (1.1) in \cite{IW2}.

Set $V=H^1(\Bbb R^n)$ and $H=L^2(\Bbb R^n)$.
Applying Theorem 1, we know that given $u^0\in V$ and
$u^1\in H$, there exists a unique $u$ satisfying
$$
\displaystyle
u\in L^2(0,\,\,T;V),\,\,
u'=\frac{du}{dt}\in L^2(0,\,\,T;V),\,\,
\frac{d}{dt}(\mbox{\boldmath $C$}(u'(\,\cdot\,))\in L^2(0,\,\,T;V' ),
$$
such that, for all $\phi\in V$
$$\displaystyle
<\frac{d}{dt}\mbox{\boldmath $C$}(u'(t)),\phi>+b_0(u'(t),v)
+a(u(t),\phi)=0,\,\,\text{a.e.}\,t\in]0,\,\,T[,
\tag {2.1}
$$
and $u(0)=u^0$ and $u'(0)=u^1$, where
$$\begin{array}{c}
\displaystyle
a(u,v)=\int_{\Bbb R^n}\nabla u\cdot\nabla vdx,\,\,b_0(u,v)=\int_{\Bbb R^n}quvdx,\,\,u,v\in V,
\end{array}
$$
and $\mbox{\boldmath $C$}:H\longrightarrow H$ is the bounded linear operator defined
by
$$\displaystyle
\mbox{\boldmath $C$}(u)=\alpha u,\,\,u\in H.
$$
Note that this $\mbox{\boldmath $C$}$ satisfies (5.11) on p. 553 in \cite{DL}
under the condition $\text{ess.inf}_{x\in\Bbb R^n}\alpha(x)>0$.
Since $q\ge 0$, $b_0$ satisfies (5.8) on p. 553 with $\beta_0=0$ in their notation.
However, equation (2.1) is {\it homogeneous}, i.e., the source term
is $0$, and by virtue of this, their proof also covers this case.

In this paper, we say that this $u$ for $u^0=0$ and
$u^1=f$ with $q=0$ and $\alpha=1$ is the weak solutions of (1.1) and (1.7), respectively.

We see that $w$ given by
$$\displaystyle
w=\int_0^Te^{-\tau t}udt
$$
belongs to $V$ and
applying integration by parts to (2.1) multiplied by $e^{-\tau T}$,
we obtain, for all $\phi\in V$
$$\displaystyle
\int_{\Bbb R^n}\{\nabla w\cdot\nabla\phi+(\tau^2\alpha+\tau q)w\phi\}dx
-\int_{\Bbb R^n}\alpha f\phi dx+e^{-\tau T}\int_{\Bbb R^3}\mbox{$\cal F$}\phi dx=0,
$$
where
$$\displaystyle
\mbox{$\cal F$}(x)=\alpha(x)u'(x,T)+(\alpha(x)\tau+q(x)) u(x,T).
$$
This means that $w$ is the weak solution of the following equation:
$$\displaystyle
(\triangle-\alpha\tau^2-q\tau)w+\alpha f
=e^{-\tau T}\mbox{$\cal F$}\,\,\text{in}\,\Bbb R^n.
\tag {2.2}
$$

\subsection{Two basic identities}

In this subsection we consider only the case when $q(x)=0$ a.e. $x\in\Bbb R^n$.
Then, it follows from (2.2) that $w$ satisfies
$$
\displaystyle
\triangle w-\alpha\tau^2 w+\alpha f=\alpha e^{-\tau T}F\,\,\,\text{in}\,\Bbb R^n,
\tag {2.3}
$$
where $\displaystyle F(x,\tau)=u'(x,T)+\tau u(x,T)$.

In what follows $f$ just belongs to $L^2(\Bbb R^n)$; $v\in H^1(\Bbb R^n)$ satisfies (1.4).

\proclaim{\noindent Proposition 2.1.}
We have
$$\begin{array}{c}
\displaystyle
\int_{\Bbb R^n}f\{(\alpha_0-\alpha)v+\alpha R\}dx
=
\tau^2\int_{\Bbb R^n}(\alpha_0-\alpha)v^2dx\\
\\
\displaystyle
+\int_{\Bbb R^n}(\vert\nabla R\vert^2+\alpha\tau^2 R^2)dx
+e^{-\tau T}\left(\int_{\Bbb R^n}\alpha FRdx-\int_{\Bbb R^n}\alpha Fvdx\right),
\end{array}
\tag {2.4}
$$
where $R=w-v$.

\endproclaim

{\it\noindent Proof.}
From (1.4) we have
$$\begin{array}{c}
\displaystyle
\int_{\Bbb R^n}\alpha_0 f wdx
=\tau^2\int_{\Bbb R^n}\alpha_0 vwdx+\int_{\Bbb R^n}\nabla v\cdot\nabla w dx.
\end{array}
$$
On the other hand, from (2.3) we have
$$\begin{array}{c}
\displaystyle
\int_{\Bbb R^n}\alpha fvdx
=\tau^2\int_{\Bbb R^n}\alpha wvdx+\int_{\Bbb R^n}\nabla w\cdot\nabla vdx
+e^{-\tau T}\int_{\Bbb R^n}\alpha Fvdx.
\end{array}
$$
Therefore we obtain
$$\begin{array}{c}
\displaystyle
\int_{\Bbb R^n}f(\alpha_0 w-\alpha v)dx
=\tau^2\int_{\Bbb R^n}(\alpha_0-\alpha)vwdx-e^{-\tau T}\int_{\Bbb R^n}\alpha Fvdx.
\end{array}
\tag {2.5}
$$
Write
$$\displaystyle
\tau^2\int_{\Bbb R^n}(\alpha_0-\alpha)vwdx
=\tau^2\int_{\Bbb R^n}(\alpha_0-\alpha)v^2dx
+\tau^2\int_{\Bbb R^n}(\alpha_0-\alpha)vRdx.
\tag {2.6}
$$
It follows from (2.3) and (1.4) that $R$ satisfies
$$\displaystyle
\triangle R-\alpha\tau^2R+(\alpha_0-\alpha)\tau^2v
-(\alpha_0-\alpha)f=\alpha e^{-\tau T}F\,\,\text{in}\,\Bbb R^n.
\tag {2.7}
$$
Thus we have
$$\begin{array}{c}
\displaystyle
\tau^2\int_{\Bbb R^n}(\alpha_0-\alpha)vRdx
=\int_{\Bbb R^n}(\vert\nabla R\vert^2+\alpha\tau^2 R^2)dx
+\int_{\Bbb R^n}(\alpha_0-\alpha)fRdx
+e^{-\tau T}\int_{\Bbb R^n}\alpha FRdx.
\end{array}
$$
A combination of this and (2.6) in (2.5) yields
$$\begin{array}{c}
\displaystyle
\int_{\Bbb R^n}f(\alpha_0 w-\alpha v)dx
=
\tau^2\int_{\Bbb R^n}(\alpha_0-\alpha)v^2dx
+\int_{\Bbb R^n}(\vert\nabla R\vert^2+\alpha\tau^2 R^2)dx\\
\\
\displaystyle
+\int_{\Bbb R^n}(\alpha_0-\alpha)fRdx
+e^{-\tau T}\left(\int_{\Bbb R^n}\alpha FRdx-\int_{\Bbb R^n}\alpha Fvdx\right).
\end{array}
$$
Since $(\alpha_0 w-\alpha v)-(\alpha_0-\alpha)R=(\alpha_0-\alpha)v+\alpha R$, we obtain (2.4).

\noindent
$\Box$

\proclaim{\noindent Proposition 2.2.}
We have
$$\begin{array}{c}
\displaystyle
\int_{\Bbb R^n}f\{(\alpha-\alpha_0)w-\alpha_0R\}dx
=
\tau^2\int_{\Bbb R^n}\frac{\alpha_0}{\alpha}(\alpha-\alpha_0)v^2dx\\
\\
\displaystyle
+\int_{\Bbb R^n}\left(\vert\nabla R\vert^2+\alpha\tau^2
\left\vert R+\left(1-\frac{\alpha_0}{\alpha}\right)v\right\vert^2\right)dx
+e^{-\tau T}\left(\int_{\Bbb R^n}\alpha FRdx+\int_{\Bbb R^n}\alpha Fvdx\right),
\end{array}
\tag {2.8}
$$
where $R=w-v$.

\endproclaim

{\it\noindent Proof.}
We recall that we have equation (2.5).
Instead of (2.6) we write
$$\displaystyle
\tau^2\int_{\Bbb R^n}(\alpha-\alpha_0)vwdx
=\tau^2\int_{\Bbb R^n}(\alpha-\alpha_0)w^2dx
+\tau^2\int_{\Bbb R^n}(\alpha_0-\alpha)wRdx.
\tag {2.9}
$$
From (2.7) we see that $R$ satisfies
$$\displaystyle
\triangle R-\alpha_0\tau^2R+(\alpha_0-\alpha)\tau^2w
-(\alpha_0-\alpha)f=\alpha e^{-\tau T}F\,\,\text{in}\,\Bbb R^n.
$$
Thus we have
$$\begin{array}{c}
\displaystyle
\tau^2\int_{\Bbb R^n}(\alpha_0-\alpha)wRdx
=\int_{\Bbb R^n}(\vert\nabla R\vert^2+\alpha_0\tau^2 R^2)dx
+\int_{\Bbb R^n}(\alpha_0-\alpha)fRdx
+e^{-\tau T}\int_{\Bbb R^n}\alpha FRdx.
\end{array}
$$
Now a combination of this and (2.9) in (2.5) 
yields
$$\begin{array}{c}
\displaystyle
\int_{\Bbb R^n}f(\alpha v-\alpha_0 w)dx
=
\tau^2\int_{\Bbb R^n}(\alpha-\alpha_0)w^2dx
+\int_{\Bbb R^n}(\nabla R\vert^2+\alpha_0\tau^2 R^2)dx\\
\\
\displaystyle
+\int_{\Bbb R^n}(\alpha_0-\alpha)fRdx
+e^{-\tau T}\left(\int_{\Bbb R^n}\alpha FRdx+\int_{\Bbb R^n}\alpha Fvdx\right).
\end{array}
\tag {2.10}
$$
Finally write $(\alpha v-\alpha_0 w)-(\alpha_0-\alpha)R=(\alpha-\alpha_0)w-\alpha_0R$ and
$$\displaystyle
\alpha_0R^2+(\alpha-\alpha_0)w^2
=\alpha\left\vert w-\frac{\alpha_0}{\alpha}v\right\vert^2+\frac{\alpha_0}{\alpha}(\alpha-\alpha_0)v^2.
$$
Combining these with (2.10) we obtain (2.8).

\noindent
$\Box$

\section{Proof of Theorem 1.1.}

First we derive two asymptotic estimates.

\proclaim{\noindent Proposition 3.1.}
We have, as $\tau\longrightarrow\infty$
$$\begin{array}{c}
\displaystyle
I_f(\tau,T)
\le
\tau^2\int_{\Bbb R^n}\frac{\alpha_0}{\alpha}(\alpha_0-\alpha)v^2dx+O(\tau^{-1}e^{-\tau T})
\end{array}
\tag {3.1}
$$
and
$$\displaystyle
I_f(\tau,T)
\ge \tau^2\int_{\Bbb R^n}(\alpha_0-\alpha)v^2dx+O(\tau^{-1}e^{-\tau T}).
\tag {3.2}
$$

\endproclaim

{\it\noindent Proof.}
We start with describing two simple estimates for $v$ and $R=w-v$.
It follows from (1.4) that
$$\begin{array}{c}
\displaystyle
\int_{\Bbb R^n}(\vert\nabla v\vert^2+\alpha_0\tau^2 v^2-\alpha_0 fv)dx=0,
\end{array}
$$
that is,
$$\begin{array}{c}
\displaystyle
\int_{\Bbb R^n}\left\{\vert\nabla v\vert^2+\alpha_0\left(\tau v-\frac{f}{2\tau}\right)^2\right\}dx
=\frac{1}{4\tau^2}\int_{\Bbb R^n}\alpha_0 f^2dx.
\end{array}
$$
It is easy to see that from this we obtain, as $\tau\longrightarrow\infty$
$$\displaystyle
\int_{\Bbb R^n}(\vert\nabla v\vert^2+\alpha_0\tau^2v^2)dx=O(\tau^{-2})
$$
and in particular,
$$
\displaystyle \Vert v\Vert_{L^2(\Bbb R^n)}=O(\tau^{-2}).
\tag {3.3}
$$

Next rewrite (2.4) as
$$\begin{array}{c}
\displaystyle
\tau^2\int_{\Bbb R^n}(\alpha_0-\alpha)v^2dx
+\int_{\Bbb R^n}
\left\{\vert\nabla R\vert^2+\alpha\left(\tau R-\frac{f-e^{-\tau t}F}{2\tau}\right)^2\right\}dx\\
\\
\displaystyle
=\int_{\Bbb R^n}(\alpha_0-\alpha)fvdx
+\frac{1}{4\tau^2}\int_{\Bbb R^n}\alpha(f-e^{-\tau T}F)^2dx.
\end{array}
$$
Since $F=u'(x,T)+\tau u(x,T)$ and $\Vert u'(\,\cdot\,,T)\Vert_{L^2(\Bbb R^n)}
+\Vert u(\,\cdot\,,T)\Vert_{L^2(\Bbb R^n)}<\infty$, it follows from this and (3.3) that
$$\displaystyle
\int_{\Bbb R^n}
\left\{\vert\nabla R\vert^2+\alpha\left(\tau R-\frac{f-e^{-\tau t}F}{2\tau}\right)^2\right\}dx
=O(\tau^{-2})
$$
and hence
$$\displaystyle
\int_{\Bbb R^n}(\vert\nabla R\vert^2+\alpha\tau^2 R^2)dx=O(\tau^{-2}).
$$
In particular, we have
$$\displaystyle
\Vert R\Vert_{L^2(\Bbb R^n)}=O(\tau^{-2}).
\tag {3.4}
$$
Applying (3.3) and (3.4) to (2.4) and (2.8), we obtain, as $\tau\longrightarrow\infty$
$$\begin{array}{c}
\displaystyle
\int_{\Bbb R^n}f\{(\alpha_0-\alpha)v+\alpha R\}dx
=
\tau^2\int_{\Bbb R^n}(\alpha_0-\alpha)v^2dx\\
\\
\displaystyle
+\int_{\Bbb R^n}(\vert\nabla R\vert^2+\alpha\tau^2 R^2)dx+O(\tau^{-1}e^{-\tau T})
\end{array}
\tag {3.5}
$$
and
$$\begin{array}{c}
\displaystyle
\int_{\Bbb R^n}f\{(\alpha-\alpha_0)w-\alpha_0R\}dx
=
\tau^2\int_{\Bbb R^n}\frac{\alpha_0}{\alpha}(\alpha-\alpha_0)v^2dx\\
\\
\displaystyle
+\int_{\Bbb R^n}\left(\vert\nabla R\vert^2+\alpha\tau^2
\left\vert R+\left(1-\frac{\alpha_0}{\alpha}\right)v\right\vert^2\right)dx
+O(\tau^{-1}e^{-\tau T}).
\end{array}
\tag {3.6}
$$
Here note that
$$\displaystyle
(\alpha_0-\alpha)
+\frac{\alpha_0}{\alpha}(\alpha-\alpha_0)
=-\frac{(\alpha-\alpha_0)^2}{\alpha}
\tag {3.7}
$$
and $\displaystyle
(\alpha_0-\alpha)v+\alpha R
+(\alpha-\alpha_0)w-\alpha_0R
=2(\alpha-\alpha_0)R$.
Thus, summing (3.5) and (3.6) up, we obtain
$$\begin{array}{l}
\displaystyle
2\int_{\Bbb R^n}f(\alpha-\alpha_0)Rdx+\tau^2\int_{\Bbb R^n}\frac{(\alpha-\alpha_0)^2}{\alpha}v^2dx
=\int_{\Bbb R^n}(\vert\nabla R\vert^2+\alpha\tau^2 R^2)dx\\
\\
\displaystyle
+\int_{\Bbb R^n}\left(\vert\nabla R\vert^2+\alpha\tau^2
\left\vert R+\left(1-\frac{\alpha_0}{\alpha}\right)v\right\vert^2\right)dx
+O(\tau^{-1}e^{-\tau T}).
\end{array}
\tag {3.8}
$$
From the assumption on $f$ we have $\text{supp}\, f=\overline B$.
Since $\alpha(x)=\alpha_0(x)$ a.e. $x\in\Bbb R^n\setminus D$ and $\overline B\cap\overline D=\emptyset$,
we have $\alpha(x)=\alpha_0(x)$ a.e. $x\in B$.
Therefore the first integral in the left-hand side on equation (3.8) vanishes.
Then (3.8) gives
$$\begin{array}{c}
\displaystyle
\int_{\Bbb R^n}(\vert\nabla R\vert^2+\alpha\tau^2 R^2)dx
+\int_{\Bbb R^n}\left(\vert\nabla R\vert^2+\alpha\tau^2
\left\vert R+\left(1-\frac{\alpha_0}{\alpha}\right)v\right\vert^2\right)dx\\
\\
\displaystyle
=\tau^2\int_{\Bbb R^n}\frac{(\alpha-\alpha_0)^2}{\alpha}v^2dx+O(\tau^{-1}e^{-\tau T})
\end{array}
$$
and hence
$$\displaystyle
\int_{\Bbb R^n}(\vert\nabla R\vert^2+\alpha\tau^2 R^2)dx
\le
\tau^2\int_{\Bbb R^n}\frac{(\alpha-\alpha_0)^2}{\alpha}v^2dx+O(\tau^{-1}e^{-\tau T}).
\tag {3.9}
$$
Now it follows from (3.5), (3.7) and (3.9) that (3.1) and (3.2) are valid.

\noindent
$\Box$

{\bf\noindent Remark 3.1.}
We note that (3.7) gives
$$\displaystyle
\alpha_0-\alpha\le\frac{\alpha_0}{\alpha}(\alpha_0-\alpha).
$$
Thus (3.1) and (3.2) are reasonable.

Now assume that $D=\emptyset$.  Then $\alpha_0=\alpha$ and thus (3.1) and (3.2) yield
$$\displaystyle
I_f(\tau,T)=O(\tau^{-1}e^{-\tau T}).
$$
This gives (i) of Theorem 1.1.

The proof of (ii) is as follows.
Since $h$ satisfies (A.I), it follows from (3.1) and (3.2) that
$$\displaystyle
I_f(\tau,T)
\le -A\tau^2\Vert v\Vert_{L^2(D)}^2
+O(\tau^{-1}e^{-\tau T}),
\tag {3.10}
$$
where $A=-(m_0^2)/(M^2)C$, $M=\sqrt{\text{ess.sup}_{x\in\Bbb R^n}\alpha(x)}$ and $C$ comes from (A.I);
$$\displaystyle
I_f(\tau,T)
\ge -A'\tau^2\Vert v\Vert_{L^2(D)}^2
+O(\tau^{-1}e^{-\tau T}),
\tag {3.11}
$$
where $A'=\Vert h\Vert_{L^{\infty}(D)}>0$.

Thus it suffices to give a lower and upper estimate of $v$ over $D$.
Given $\lambda>0$
define
$$
\displaystyle
G_{\lambda}(\xi)=\left\{
\begin{array}{ll}
\displaystyle
\frac{1}{2\lambda}e^{-\lambda\vert\xi\vert}, & \mbox{if $n=1$,}\\
\\
\displaystyle
\frac{1}{2\pi}K_0(\lambda\vert \xi\vert), & \mbox{if $n=2$,}\\
\\
\displaystyle
\frac{e^{-\lambda\vert\xi\vert}}{4\pi\vert\xi\vert}, & \mbox{if $n=3$},
\end{array}
\right.
\tag {3.12}
$$
where $K_0$ is the modified Bessel function of the ssecond kind of order $0$ (see \cite{O}).
It seems that the following lemma is closely related to the
maximum principle or comparison principle \cite{GT}. However, our
final purpose is to consider the electromagnetic wave which
satisfies a {\it system}. So in Appendix we give a proof without
making use of such principles.

\proclaim{\noindent Lemma 3.1.}
Let $f\in L^2(\Bbb R^n)$ and satisfy
$f(x)\ge 0$ a.e. $x\in\Bbb R^n$.
Let $v\in H^1(\Bbb R^n)$ be the weak solution of (1.4).
We have
$$\displaystyle
v(x)\ge\int_{\Bbb R^n}\alpha_0(y)f(y)G_{M_0\tau}(x-y)dy\,\,\text{a.e.}\,x\in\Bbb R^n,
\tag {3.13}
$$
and
$$
\displaystyle
v(x)\le\int_{\Bbb R^n}\alpha_0(y)f(y)G_{m_0\tau}(x-y)dy\,\,\text{a.e.}\,x\in\Bbb R^n.
\tag {3.14}
$$

\endproclaim

Let us continue the proof of (ii).
Since the case when $n=1, 2$ can be treated easily, hereafter
we only consider the case when $n=3$.
By the mean value theorem \cite{CH} we have
$$\displaystyle
\frac{1}{4\pi}
\int_B\frac{e^{-\tau M_0\vert x-y\vert}}{\vert x-y\vert}dy
=\frac{\varphi(\tau M_0\eta)}{(M_0\tau)^3}
\frac{e^{-\tau M_0\vert x-p\vert}}{\vert x-p\vert},\,\,x\in\Bbb R^3\setminus\overline B,
\tag {3.15}
$$
where $p$ and $\eta$ are the center and radius of $B$,
respectively and $\varphi(\xi)=\xi\cosh\xi-\sinh\xi$.

By \cite{IkIt2}, we know that there exists a positive constant $C'$ and number $\mu\in\Bbb R$ such that, for all $\tau>>1$
$$\displaystyle
\tau^{\mu}e^{2\tau d_{\partial D}(p)}
\int_D\frac{e^{-2\tau\vert x-p\vert}}{\vert x-p\vert^2}dx
\ge C'.
$$
Applying this and (3.15) to a lower bound derived from the lower bounds for $\alpha_0$ and $f$ on $B$
for the right-hand side on (3.13), we obtain,
for all $\tau>>1$
$$\begin{array}{c}
\displaystyle
\Vert v\Vert_{L^2(D)}^2
\ge
C^2\frac{m^4\varphi(M_0\tau\eta)^2}{(M_0\tau)^{6+\mu}}e^{-2M_0\tau d_{\partial D}(p)}
\times
(M_0\tau)^{\mu}e^{2M_0\tau d_{\partial D}(p)}
\int_D\frac{e^{-2M_0\tau\vert x-p\vert}}{\vert x-p\vert^2}dx\\
\\
\displaystyle
\ge
C^2C'\frac{m^4\varphi(M_0\tau\eta)^2}{(M_0\tau)^{6+\mu}}e^{-2M_0\tau d_{\partial D}(p)},
\end{array}
\tag {3.16}
$$
where $C=\text{ess.inf}_{x\in B}\,f(x)$.
Since as $\xi\longrightarrow\infty$
$\displaystyle\varphi(\xi)\sim e^{\xi}/2$,
from (3.16) we obtain, for all $\tau>>1$
$\displaystyle
\Vert v\Vert_{L^2(D)}^2
\ge C''\tau^{-(4+\mu)}e^{-2M_0\tau(d_{\partial D}(p)-\eta)}$,
where $C''$ is a positive constant.  Since $d_{\partial D}(p)-\eta=\text{dist}\,(D,B)$, we finally obtain,
for all $\tau>>1$
$$\displaystyle
\Vert v\Vert_{L^2(D)}^2
\ge C''\tau^{-(4+\mu)}e^{-2M_0\tau\text{dist}\,(D,B)}.
\tag {3.17}
$$

Now a combination of this and (3.10) yields $\displaystyle\lim_{\tau\longrightarrow\infty}e^{\tau T}I_f(\tau,T)=-\infty$
if $T>2M_0\text{dist}\,(D,B)$.

Moreover, it is easy to obtain from (3.14) that, for a positive constant $C^{'''}$ and
all $\tau>0$ we have
$$\displaystyle
\Vert v\Vert_{L^2(D)}^2
\le C^{'''}e^{-2m_0\tau\text{dist}\,(D,B)}.
\tag {3.18}
$$
Now it is easy to verify that a combination of (3.10) and (3.17) yields (1.5);
a combination of (3.11) and (3.18) yields (1.6).

A similar argument based on (3.1) and (3.2) works also for the case when $h$ satisfies (A.II).

\section{Proof of Theorem 1.2}

From (2.2) we have
$$\displaystyle
(\triangle-\tau^2-\tau q)w+f=e^{-\tau T}\mbox{$\cal F$}\,\,\text{in}\,\Bbb R^n
\tag {4.1}
$$
and $\mbox{$\cal F$}(x,\tau)=u'(x,T)+(\tau+q(x))u(x,T)$.

Define
$$
\displaystyle
\tilde{\alpha}=\tilde{\alpha}(x,\tau)=1+\frac{q}{\tau},\,\,
\tilde{f}=\tilde{f}(x,\tau)
=\left(1+\frac{q}{\tau}\right)^{-1}f,\,\,
\tilde{F}=\tilde{F}(x,\tau)
=\left(1+\frac{q}{\tau}\right)^{-1}\mbox{$\cal F$}.
$$
Then (4.1) becomes
$$\displaystyle
(\triangle-\tau^2\tilde{\alpha})w+\tilde{\alpha}\tilde{f}
=\tilde{\alpha}e^{-\tau T}\tilde{F}\,\,\text{in}\,\Bbb R^n.
$$

Define
$$\displaystyle
\displaystyle
\tilde{\alpha_0}=\tilde{\alpha_0}(x,\tau)
=1+\frac{q_0}{\tau},\,\,
\tilde{f_0}=\tilde{f_0}(x,\tau)
=\left(1+\frac{q_0}{\tau}\right)^{-1}f.
$$
Let $v\in H^1(\Bbb R^n)$ be the solution of (1.9).
Then (1.9) becomes
$$\displaystyle
(\triangle-\tilde{\alpha_0}\tau^2)v+\tilde{\alpha_0}\tilde{f_0}=0\,\,\text{in}\,\Bbb R^n.
\tag {4.2}
$$
Note that $\tilde{\alpha_0}\tilde{f_0}=f=\tilde{\alpha}\tilde{f}$ and $\tilde{\alpha}\tilde{F}=\mbox{$\cal F$}$.

Define $R=w-v$.
We start with the following two integral identities which can be derived along the same line
as that of Propositions 2.2 and 2.3 and thus we omit the description of their proofs:
$$\begin{array}{c}
\displaystyle
J_f(\tau,T)
=\tau^2\int_{\Bbb R^n}(\tilde{\alpha_0}-\tilde{\alpha})v^2dx
+\int_{\Bbb R^n}(\vert\nabla R\vert^2+\tilde{\alpha}\tau^2 R^2)dx\\
\\
\displaystyle
+e^{-\tau T}
\left(\int_{\Bbb R^n}\mbox{$\cal F$}Rdx
-\int_{\Bbb R^n}\mbox{$\cal F$}vdx\right);
\end{array}
\tag {4.3}
$$
$$\begin{array}{c}
\displaystyle
-J_f(\tau,T)
=\tau^2\int_{\Bbb R^n}\frac{\tilde{\alpha_0}}{\tilde{\alpha}}(\tilde{\alpha}-\tilde{\alpha_0})v^2dx
\\
\\
\displaystyle
+\int_{\Bbb R^n}\left(\vert\nabla R\vert^2+\tilde{\alpha}\tau^2\left\vert R+\left(1-\frac{\tilde{\alpha_0}}{\tilde{\alpha}}\right)v\right\vert^2\right)dx
+e^{-\tau T}
\left(\int_{\Bbb R^n}\mbox{$\cal F$}Rdx+\int_{\Bbb R^n}\mbox{$\cal F$}vdx\right).
\end{array}
\tag {4.4}
$$
Note also that
$$\displaystyle
\tau^2\int_{\Bbb R^n}(\tilde{\alpha_0}-\tilde{\alpha})v^2dx
=\tau\int_{\Bbb R^n}(q_0-q)v^2dx
$$
and
$$\displaystyle
-\tau^2\int_{\Bbb R^n}\frac{\tilde{\alpha_0}}{\tilde{\alpha}}(\tilde{\alpha}-\tilde{\alpha_0})v^2dx
=\tau\int_{\Bbb R^n}\frac{\tau+q_0}{\tau+q}(q_0-q)v^2dx.
$$
Then, applying a similar argument as done in the proof of Proposition 3.1 to (4.3) and (4.4) we obtain

\proclaim{\noindent Proposition 4.1.}
We have, as $\tau\longrightarrow\infty$
$$\begin{array}{c}
\displaystyle
J_f(\tau,T)
\le
\tau\int_{\Bbb R^n}\frac{\tau+q_0}{\tau+q}(q_0-q)v^2dx+O(\tau^{-1}e^{-\tau T})
\end{array}
\tag {4.5}
$$
and
$$\displaystyle
J_f(\tau,T)
\ge\tau\int_{\Bbb R^n}(q_0-q)v^2dx+O(\tau^{-1}e^{-\tau T}).
\tag {4.6}
$$

\endproclaim

Thus it suffices to prepare the following lower and upper estimates for $v$.

\proclaim{\noindent Lemma 4.1.}
Let $f\in L^2(\Bbb R^n)$ and satisfy
$f(x)\ge 0$ a.e. $x\in\Bbb R^n$.
Let $v\in H^1(\Bbb R^n)$ be the weak solution of (1.9).
We have
$$\displaystyle
v(x)\ge\int_{\Bbb R^n}f(y)G_{L_0(\tau)\tau}(x-y)dy\,\,\text{a.e.}\,x\in\Bbb R^n
\tag {4.7}
$$
and
$$\displaystyle
v(x)\le \int_{\Bbb R^n}f(y)G_{\tau}(x-y)dy\,\,\text{a.e.}\,x\in\Bbb R^n,
\tag {4.8}
$$
where $G_{L_0(\tau)\tau}=G_{\lambda}\vert_{\lambda=L_0(\tau)\tau}$, $G_{\tau}=G_{\lambda}\vert_{\lambda=\tau}$
with $G_{\lambda}$ given by (3.12),
$$\displaystyle
L_0(\tau)=\sqrt{1+\frac{L_0}{\tau}}
\tag {4.9}
$$
and $L_0=\text{ess.sup}_{x\in\Bbb R^n}\,q_0(x)$.
\endproclaim

{\it\noindent Proof.}
Since $v$ satisfies (4.2) and $L_0(\tau)$ given by (4.9) satisfies $\tilde{\alpha_0}(x)\le L_0(\tau)^2$ a.e. $x\in\Bbb R^n$,
from (3.13) in Lemma 3.1 one obtains (4.7).
The proof of (4.8) is as follows.
$v$ has the expression
$\displaystyle
v=v_0-\epsilon_0$,
where $v_0\in H^1(\Bbb R^n)$ solves
$\displaystyle
(\triangle-\tau^2)v_0+f=0\,\,\text{in}\,\Bbb R^n$
and $\epsilon_0\in H^1(\Bbb R^n)$ solves
$$\displaystyle
(\triangle-\tau^2-\tau q_0)\epsilon_0+\tau q_0v_0=0\,\,\text{in}\,\Bbb R^n.
\tag {4.10}
$$
$v_0$ has the explicit form
$$\displaystyle
v_0(x)=\int_{\Bbb R^n}f(y)G_{\tau}(x-y)dy\ge 0
$$
and hence $\tau q_0v_0\ge 0$.
Applying (4.7) to (4.10) we obtain
$$\displaystyle
\epsilon_0(x)\ge\int_{\Bbb R^n}\tau q_0(y)v_0(y)G_{L_0(\tau)\tau}(x-y)dy
$$
and hence $\epsilon_0(x)\ge 0$ a.e. $x\in\Bbb R^n$.  Therefore we obtain $v\le v_0$ and thus (4.8).

\noindent
$\Box$

Now noting $L_0(\tau)\longrightarrow 1$ as $\tau\longrightarrow\infty$, it is not difficult to deduce all the conclusions in Theorem 1.2
from (4.5), (4.6) and Lemma 4.1 as done in the proof of Theorem 1.1.

\section{Remarks and further problems}

\subsection{How to compute the indicator function without knowledge about $\alpha_0$ outside $B$}

In this work we think that $D$ is an invader into a space with the refractive index $\alpha_0$.
However, it will be difficult to know the detail of $\alpha_0$ and find $v$ on $B$ which is the solution
of equation (1.4).  In this section we describe an experimental computation
procedure of the indicator function from the observed data in the two spaces one of which has an invader
and another does not have an invader yet.

Let $V=V(x,t)$ be the weak solution of
$$\displaystyle
\left\{
\begin{array}{ll}
\displaystyle
\alpha_0(x)\partial_t^2V-\triangle V=0 & \text{in}\,\Bbb R^n\times\,]0,\,T[,
\\
\\
\displaystyle
V(x,0)=0 & \text{in}\,\Bbb R^n,\\
\\
\displaystyle
\partial_tV(x,0)=\chi_B(x) & \text{in}\,\Bbb R^n.
\end{array}
\right.
\tag {5.1}
$$
Generate wave $V_e$ by the initial data $f=\chi_B$ in the space which has no invader
and observe $V_e$ on $B$ over time interval $]0,\,T[$.
We assume that $V_e$ on $B$ is given by $V$ on $B$.

From Proposition 3.1 in the case when $\alpha_0=\alpha$ we obtain
$$\displaystyle
\int_B\alpha_0(v_e-v)dx=O(\tau^{-1}e^{-\tau T}),
$$
where
$$\displaystyle
v_e(x,\tau)=\int_0^Te^{-\tau t}V_e(x,t)dt.
$$
Thus, we obtain
$$\displaystyle
I_f(\tau,T)_e\equiv
\int_B\alpha_0(w-v_e)dx
=I_f(\tau,T)+O(\tau^{-1}e^{-\tau T}).
$$
Then, it follows from Proposition 3.1 that, as $\tau\longrightarrow\infty$
$$\begin{array}{c}
\displaystyle
I_f(\tau,T)_e
\le
\tau^2\int_{\Bbb R^n}\frac{\alpha_0}{\alpha}(\alpha_0-\alpha)v^2dx+O(\tau^{-1}e^{-\tau T})
\end{array}
$$
and
$$\displaystyle
I_f(\tau,T)_e
\ge \tau^2\int_{\Bbb R^n}(\alpha_0-\alpha)v^2dx+O(\tau^{-1}e^{-\tau T}).
$$

Therefore, one can transplant all the results in Theorem 1.1 into the present case
and we obtain

\proclaim{\noindent Theorem 5.1.}
We have:

(i) if $D=\emptyset$, then for all $T>0$ it holds that
$$\displaystyle
\lim_{\tau\longrightarrow\infty}e^{\tau T}I_f(\tau,T)_e=0;
$$

(ii) if $D\not=\emptyset$ and $h$ satisfies (A.I), then for all $T>2M_0\text{dist}(D,B)$ it holds that
$$\displaystyle
\lim_{\tau\longrightarrow\infty}e^{\tau T}I_f(\tau,T)_e=-\infty;
$$

(iii)  if $D\not=\emptyset$ and $h$ satisfies (A.II), then for all $T>2M_0\text{dist}(D,B)$ it holds that
$$\displaystyle
\lim_{\tau\longrightarrow\infty}e^{\tau T}I_f(\tau,T)_e=\infty.
$$

Moreover, in case of both (ii) and (iii) we have, for all $T>2M_0\text{dist}(D,B)$
$$\left\{
\displaystyle
\begin{array}{l}
\displaystyle
\liminf_{\tau\longrightarrow\infty}
\frac{1}{2\tau}
\log\left\vert I_f(\tau,T)_e\right\vert\ge
-M_0\text{dist}\,(D,B),\\
\\
\displaystyle
\limsup_{\tau\longrightarrow\infty}
\frac{1}{2\tau}
\log\left\vert I_f(\tau,T)_e\right\vert
\le -m_0\text{dist}\,(D,B).
\end{array}
\right.
$$

\endproclaim

Note that, if we know $\alpha_0$ on $B$, then one can compute the indicator function
$\tau\longmapsto I_f(\tau,T)_e$
from the {\it experimental data} and the values of $\alpha_0$ on $B$.

We just need the following qualitative knowledge:

(i)  the governing equation of the observed wave in the space which has a penetrable obstacle
takes the form (1.1) and its refractive index is given by (1.2);

(ii) the governing equation of the observed wave in the space which has no obstacle yet takes the form
(5.1) with $\alpha_0$ in (1.2).

Summing up, we can say that: one can know the existence of something added to a reference space
by comparing the ``snap shot'' $u$ on $B$ with a ``reference snap shot'' $V_e$ on $B$
even in the case when: the reference space has a complicated rough refractive index $\alpha_0$;
$\alpha_0$ is {\it unknown} outside $B$.

\subsection{One-space dimensional case}

Finally let us describe one non trivial application of the method presented in the proof of Theorem 1.1 in the case when
the space dimension is one and the background medium is not homogeneous.

We assume that $\alpha_0$ is piecewise constant and takes the form
$$\displaystyle
\alpha_0(x)=\left\{
\begin{array}{ll}
\displaystyle
1, & \mbox{if $x<a$ or $b<x$}\\
\\
\displaystyle
k_0, & \mbox{if $a<x<b$},
\end{array}
\right.
\tag {5.2}
$$
where $-\infty<a<b<\infty$ and $k_0$ is a positive constant.

We choose $f=\chi_B$, where $\chi_B$ denotes the characteristic
function of open interval $B=]p-\epsilon,\,p+\epsilon[$ with a
fixed $p$ satisfying $p+\epsilon<a$. $f$ is a simple model of the
disturbance given at $t=0$ from the left side of the wall
$[a,\,b]$.

We assume that $D=]c,\,d[$ with, for simplicity $b<c<d<\infty$.
This means that obstacle $D$ is located behind the wall $]a,\,b[$
from the observer.

Define
$$\displaystyle
\varphi
=a-(p+\epsilon)+\sqrt{k_0}(b-a)+(c-b).
$$
The quantity $2\varphi$ coincides with the time of flight of the
signal which propagates as
$$\displaystyle
x_0=p+\epsilon\longrightarrow x_1=a\longrightarrow x_2=b\longrightarrow x_3=c\longrightarrow x_4=b\longrightarrow x_5=a\longrightarrow x_6=p+\epsilon,
$$
where the propagation speed of the signal in $]p+\epsilon,\, a[$ and $]b,\,c[$ is $1$, in $]a,\,b[$ is $1/\sqrt{k_0}$.

Note that knowing $c$ is equivalent to knowing $\varphi$ provided
the wall thickness $b-a$ and the propagation speed of the wave inside the wall $1/\sqrt{k_0}$ are known.

Let $\alpha$ be the same as (1.2) with $n=1$, $h\in L^{\infty}(D)$ and $\alpha_0$ given by (5.2).
Let $v$ be the weak solution of (1.4) with $f=\chi_B$.
In this case we have the expression
$$\displaystyle
I_f(\tau,T)=\int_B(w-v)dx.
$$

\proclaim{\noindent Theorem 5.2.}
Let $T>2\varphi$.
We have:

(i)  if $h$ satisfies (A.I), then
$$\displaystyle
\lim_{\tau\longrightarrow\infty}e^{\tau T}I_f(\tau,T)=-\infty;
$$

(ii)  if $h$ satisfies (A.II), then
$$\displaystyle
\lim_{\tau\longrightarrow\infty}e^{\tau T}I_f(\tau,T)=\infty.
$$

Moreover, in case of both (i) and (ii) we have
$$\displaystyle
\lim_{\tau\longrightarrow\infty}
\frac{1}{2\tau}
\log\left\vert I_f(\tau,T)\right\vert
=-\varphi.
\tag {5.3}
$$

\endproclaim

The proof is based on the following asymptotic formula of the
solution of (1.4) which corresponds to Lemma 3.1 and the argument
of the proof done in that of Theorem 1.1 to the present case.

\proclaim{\noindent Lemma 5.1.}
We have, as $\tau\longrightarrow\infty$
$$\displaystyle
2\tau e^{2\tau\varphi}\int_D v^2dx=1+O(\tau^{-2}).
\tag {5.4}
$$
\endproclaim

{\it\noindent Proof.}
A direct computation shows that
$v$ has the expression
\[
v(x,\tau)=\left\{
\begin{array}{ll}
\displaystyle
A(\tau)e^{\tau x}, & \mbox{if $x<p-\epsilon$}\\
\\
\displaystyle
B(\tau)e^{\tau x}+C(\tau)e^{-\tau x}+\frac{1}{\tau^2}, & \mbox{if $p-\epsilon<x<p+\epsilon$}\\
\\
\displaystyle
D(\tau)e^{\tau x}+G(\tau)e^{-\tau x}, & \mbox{if $p+\epsilon<x<a$}\\
\\
\displaystyle
H(\tau)e^{\sqrt{k_0}\,\tau x}+K(\tau)e^{-\sqrt{k_0}\,\tau x}, & \mbox{if $a<x<b$}\\
\\
\displaystyle
L(\tau)e^{-\tau x}, & \mbox{if $b<x$,}
\end{array}
\right.
\]
where
$$\begin{array}{l}
\displaystyle
A(\tau)=D(\tau)+\frac{e^{-\tau(p-\epsilon)}-e^{-\tau(p+\epsilon)}}{2\tau^2},\,\,
B(\tau)=A(\tau)-\frac{e^{-\tau(p-\epsilon)}}{2\tau^2},
\\
\\
\displaystyle
C(\tau)=-\frac{e^{\tau(p-\epsilon)}}{2\tau^2},\,\,
G(\tau)
=\frac{e^{\tau(p+\epsilon)}-e^{\tau(p-\epsilon)}}{2\tau^2};\\
\\
\displaystyle
H(\tau)=\frac{\sqrt{k_0}+1}{2\sqrt{k_0}}D(\tau)e^{\tau(1-\sqrt{k_0})a}
+\frac{\sqrt{k_0}-1}{2\sqrt{k_0}}G(\tau)e^{-\tau(1+\sqrt{k_0})a};\\
\\
\displaystyle
K(\tau)=\frac{\sqrt{k_0}-1}{2\sqrt{k_0}}D(\tau)e^{\tau(1+\sqrt{k_0})a}
+\frac{\sqrt{k_0}+1}{2\sqrt{k_0}}G(\tau)e^{-\tau(1-\sqrt{k_0})a};
\\
\\
\displaystyle
D(\tau)
=-
\frac{(k_0-1)e^{-\tau(a-p-\epsilon)}e^{-\tau a}(1+e^{-2\tau\sqrt{k_0}(b-a)})(1-e^{-2\tau\epsilon})}
{\displaystyle
2\tau^2(\sqrt{k_0}+1)^2\left\{
1-\left(\frac{\sqrt{k_0}-1}{\sqrt{k_0}+1}\right)^2e^{-2\sqrt{k_0}\tau(b-a)}\right\}};
\\
\\
\displaystyle
L(\tau)=e^{\tau b}e^{-\tau\sqrt{k_0}(b-a)}e^{-\tau(a-p-\epsilon)}\\
\\
\displaystyle
\times
\left\{
\frac{1}{1-e^{-2\tau\epsilon}}
-\frac{1}{2\tau^2}
\left(\frac{\sqrt{k_0}-1}{\sqrt{k_0}+1}\right)^2
\frac{(1+e^{-2\tau\sqrt{k_0}(b-a)})(1-e^{-2\tau\epsilon})}
{\displaystyle
1-\left(\frac{\sqrt{k_0}-1}{\sqrt{k_0}+1}\right)^2e^{-2\sqrt{k_0}\tau(b-a)}}\right\}.
\end{array}
$$
Thus, as $\tau\longrightarrow\infty$
$$\displaystyle
L(\tau)=e^{\tau b}e^{-\tau\sqrt{k_0}(b-a)}e^{-\tau(a-p-\epsilon)}
(1+O(\tau^{-2})).
\tag {5.5}
$$
Since $v(x)=L(\tau)e^{-\tau x}$ for $x>b$, we have
$$\displaystyle
\int_c^d v^2dx
=\frac{(L(\tau)e^{-\tau c})^2}{2\tau}(1-e^{-2\tau(d-c)}).
$$
Then (5.5) yields (5.4).

\noindent
$\Box$

Theorem 5.2 suggests that, if $v\in H^1(\Bbb R^n)$ is the solution
of (1.4) and one knows the leading term of $\Vert v\Vert_{L^2(D)}$
as $\tau\longrightarrow\infty$ like (5.4) in Lemma 5.1, then one
can obtain a formula in three-dimensions like (5.3) in Theorem 5.2
instead of estimates (1.5) and (1.6). For the determination of the
leading term, usually, one has to solve the {\it eikonal equation}
$\vert\nabla\Psi\vert^2=\alpha_0(x)$.  In some restricted cases it
is possible to solve the equation, see \cite{U}. The case treated
in Theorem 5.2 is just simplest one of such cases. However, we do
not go into such research direction further since the solvability
of the eikonal equation requires some regularity for $\alpha_0$
and we are seeking a method that works for finding an obstacle
embedded in a {\it rough} background medium.  We think that the
method presented in this paper is the first one for the purpose.
In a forthcoming paper we will consider the original problem which
is formulated by the Maxwell system.

$$\quad$$

\centerline{{\bf Acknowledgments}}

The author thanks the referees for their valuable and constructive comments. 
The author was partially supported by Grant-in-Aid for
Scientific Research (C)(No. 25400155) of Japan  Society for
the Promotion of Science.

$$\quad$$

\section{Appendix.  Proof of Lemma 3.1}

We make use of the following elementary fact and an iteration process.

\proclaim{\noindent Lemma A.}  Let $\lambda>0$.
Given $f\in L^2(\Bbb R^n)$ let $v\in H^1(\Bbb R^n)$ be the weak solution of
$$\displaystyle
(\triangle-\lambda^2)v+f=0\,\,\text{in}\,\Bbb R^n.
$$
Then we have
$$\displaystyle
\Vert v\Vert_{L^2(\Bbb R^n)}\le\lambda^{-2}\Vert f\Vert_{L^2(\Bbb R^n)}
\tag {A.1}
$$
and
$$\displaystyle
\Vert\nabla v\Vert_{L^2(\Bbb R^n)}\le (2\lambda)^{-1}\Vert f\Vert_{L^2(\Bbb R^n)}.
\tag {A.2}
$$

\endproclaim

{\it\noindent Proof of (3.13).}
In what follows, $G_{M_0\tau}=G_{\lambda}\vert_{\lambda=M_0\tau}$ with $G_{\lambda}$ given by (3.12).

Rewrite (1.4) as
$$\displaystyle
\{\triangle-(M_0\tau)^2\}v+\left\{\alpha_0f+\tau^2(M_0^2-\alpha_0)v\right\}=0\,\,\text{in}\,\Bbb R^n.
\tag {A.3}
$$
Let $v_1\in H^1(\Bbb R^n)$ be the weak solution of $\displaystyle
\{\triangle-(M_0\tau)^2\}v_1+\alpha_0f=0\,\,\text{in}\,\Bbb R^n$.
Since $v_1$ has the expression
$$\displaystyle
v_1(x)=\int_{\Bbb R^3}\alpha_0(y)f(y)G_{M_0\tau}(x-y)dy,
$$
we have $v_1(x)\ge 0$.

Let $j=1,\cdots$.  Given $v_j\in H^1(\Bbb R^n)$ let $v_{j+1}\in H^1(\Bbb R^n)$ be the weak solution of
$$\displaystyle
\{\triangle-(M_0\tau)^2\}v_{j+1}+\{\alpha_0f+\tau^2(M_0^2-\alpha_0)v_j\}=0\,\,\text{in}\,\Bbb R^n.
\tag {A.4}
$$
Then,  $v_{j+1}-v_j$ for $j\ge 2$ satisfies
$\displaystyle
\{\triangle-(M_0\tau)^2\}(v_{j+1}-v_j)
=-\tau^2(M_0^2-\alpha_0)(v_j-v_{j-1})\,\,\text{in}\,\Bbb R^n$
and from (A.1) we obtain
$$\displaystyle
\Vert v_{j+1}-v_j\Vert_{L^2(\Bbb R^n)}
\le \left(1-\frac{m_0^2}{M_0^2}\right)\Vert v_j-v_{j-1}\Vert_{L^2(\Bbb R^n)}
$$
and hence, for $j=1,\cdots$
$$\displaystyle
\Vert v_{j+1}-v_j\Vert_{L^2(\Bbb R^n)}
\le
\left(1-\frac{m_0^2}{M_0^2}\right)^{j-1}\Vert v_2-v_1\Vert_{L^2(\Bbb R^n)}.
$$
Similarly we have
$$\displaystyle
\Vert v_2-v_1\Vert_{L^2(\Bbb R^n)}
\le\left(1-\frac{m_0^2}{M_0^2}\right)\Vert v_1\Vert_{L^2(\Bbb R^n)}
$$
and applying (A.1) to $v_1$ on this right-hand side, we obtain
$$\displaystyle
\Vert v_2-v_1\Vert_{L^2(\Bbb R^n)}
\le
\left(1-\frac{m_0^2}{M_0^2}\right)\frac{1}{\tau^2}\Vert f\Vert_{L^2(\Bbb R^n)}.
$$
Thus we obtain, for $j=1,\cdots$
$$\displaystyle
\Vert v_{j+1}-v_j\Vert_{L^2(\Bbb R^n)}
\le
\left(1-\frac{m_0^2}{M_0^2}\right)^{j}
\frac{1}{\tau^2}\Vert f\Vert_{L^2(\Bbb R^n)}.
$$
Moreover, from (A.2) we have
$$\displaystyle
\Vert\nabla(v_{j+1}-v_j)\Vert_{L^2(\Bbb R^n)}
\le
\frac{\tau(M_0^2-m_0^2)}{2M_0}\Vert v_j-v_{j-1}\Vert_{L^2(\Bbb R^n)}
$$
and similarly
$$\displaystyle
\Vert\nabla(v_{j+1}-v_j)\Vert_{L^2(\Bbb R^n)}
\le\frac{M_0}{2\tau}\left(1-\frac{m_0^2}{M_0^2}\right)^j\Vert f\Vert_{L^2(\Bbb R^n)}.
$$
Therefore, the sequence $\{v_j\}$ in $H^1(\Bbb R^n)$ converges to
$$\displaystyle
v\equiv\sum_{j=1}^{\infty}(v_{j+1}-v_j)+v_1\,\,\text{in}\,H^1(\Bbb R^n)
\tag {A.5}
$$
and we have:
$$\displaystyle
\Vert v-v_1\Vert_{L^2(\Bbb R^n)}
\le
\frac{M_0^2}{\tau^2m_0^2}\left(1-\frac{m_0^2}{M_0^2}\right)\Vert f\Vert_{L^2(\Bbb R^n)};
$$
$$\displaystyle
\Vert\nabla v-\nabla v_1\Vert_{L^2(\Bbb R^n)}
\le\frac{M_0^3}{2\tau m_0^2}\left(1-\frac{m_0^2}{M_0^2}\right)\Vert f\Vert_{L^2(\Bbb R^n)}.
$$
Thus taking the limit of (A.4) as $j\longrightarrow\infty$, we see
that $v$ satisfies (A.3) and thus is the weak solution of (1.4).

Since $v_{2}-v_1$ has the expression
$$\displaystyle
v_{2}(x)-v_{1}(x)
=\tau^2\int_{\Bbb R^3}(M_0^2-\alpha_0(y))v_1(y)G_{M_0\tau}(x-y)dy,
$$
we have $v_2(x)-v_1(x)\ge 0$.  For $j\ge 2$ we have also
$$\displaystyle
v_{j+1}(x)-v_{j}(x)
=\tau^2\int_{\Bbb R^3}(M_0^2-\alpha_0(y))(v_j(y)-v_{j-1}(y))G_{M_0\tau}(x-y)dy
$$
and thus by induction we obtain, for all $j\ge 2$
$v_{j+1}(x)-v_j(x)\ge 0$. Therefore from (A.5) and the almost
convergence property of a subsequence of
$\{\sum_{j=1}^m(v_{j+1}-v_j)\}_{m=1}^{\infty}$ we conclude that
$\displaystyle v(x)\ge v_1(x)\,\,\text{a.e.}\,x\in\Bbb R^n$. This
completes the proof of (3.13).

\noindent
$\Box$

{\it\noindent Proof of (3.14).}
Let $v_0\in H^1(\Bbb R^n)$ solve
$$\displaystyle
\{\triangle-(m_0\tau)^2\}v_0+\alpha_0f=0\,\,\text{in}\,\Bbb R^n.
$$
$v_0$ has the expression
$$\displaystyle
v_0(x)=\int_{\Bbb R^n}\alpha_0(y)f(y)G_{m_0\tau}(x-y)dy,
$$
where $G_{m_0\tau}=G_{\lambda}\vert_{\lambda=m_0\tau}$.
Thus $v_0\ge 0$.
Then, $v$ has the expression $\displaystyle v=v_0-\epsilon_0$,
where $\epsilon_0\in H^1(\Bbb R^n)$ solves
$\displaystyle\triangle\epsilon_0-\alpha_0\tau^2\epsilon_0+\tau^2(\alpha_0-m_0^2)v_0=0\,\,\text{in}\,\Bbb R^n$.
Since $\tau^2(\alpha_0-m_0^2)v_0\ge 0$, applying (3.13) to the equation above, we obtain $\epsilon_0\ge 0$.
Therefore we obtain $v\le v_0$.

\noindent
$\Box$

\vskip1cm
\noindent
e-mail address

ikehata@amath.hiroshima-u.ac.jp

\end{document}